\documentclass{article}

\usepackage[T1]{fontenc}
\usepackage[utf8]{inputenc}

\usepackage{amsthm}
\usepackage{amssymb}
\usepackage{amsmath}

\DeclareMathOperator{\LCM}{H}

\author{Michael Hellus}

\begin{document}

\title{Generalization of a connectedness result to cohomologically complete intersections}
\date{\today}
\maketitle

\noindent\texttt{michael.hellus@mathematik.uni-regensburg.de, Fakult\"at f\"ur Mathematik, Universit\"at Regensburg, D-93040, Germany}

{\let\thefootnote\relax}
{\let\thefootnote\relax}
\begin{abstract}
It is a well-known result that, in projective space over a field, every set-theoretical complete intersection of positive dimension in connected in codimension one (Hartshorne \cite[3.4.6]{H1} or \cite[Theorem 1.3]{H2}). Another important connectedness result is that a local ring with disconnected punctured sprectrum has depth at most $1$ (\cite[Proposition 2.1]{H1}). The two results are related, Hartshorne calls the latter ``the keystone to the proof'' of the former (loc.\,cit).

In this short note we show how the latter result generalizes smoothly from \emph{set-theore\-ti\-cal} to \emph{cohomologically} complete intersections, i.\,e. to ideals for which there is in terms of local cohomology no obstruction to be a complete intersection (\cite{HeSc1}, \cite{HeSc2}).

The proof is based on the fact that, for cohomologically complete intersections over a complete local ring, the endomorphism ring of the (only) local cohomology cohomology module is the ring itself (\cite[Theorem 2.2\,(iii)]{HeSt}) and hence indecomposable as a module.

\end{abstract}

\section{Introduction}

In \cite[2.1]{H1}, Hartshorne proved (then by elementary means)

\vspace{.2cm}

\noindent\textbf{Theorem 1} Let $(R,m)$ be a noetherian local ring whose punctured spectrum $\operatorname{Spec}(R)\setminus\{m\}$ is disconnected.
\[\operatorname{depth}(R)=\operatorname{depth}(m,R)\leq1\]
\hfill$\square$

\vspace{.2cm}

\noindent Hartshorne's result can be generalized:

\vspace{.2cm}

\noindent\textbf{Theorem 2} Let $R$ be a noetherian ring, indecomposable as  a module over itself (e.\,g. local, as a ring), $I,J\subsetneq R$ ideals of $R$ such that $\sqrt{I\cap J}=\sqrt0$, but $\sqrt I,\sqrt J\neq\sqrt0$.
\[\operatorname{depth}(I+J,R)\leq1\]
\hfill$\square$

To the best of my knowledge this was first observed by Irving Kaplansky (no reference). It is nowadays well known that theorem 2 can be proved quickly using local cohomology by a straightforward adaption of the method used in \cite[proof of Proposition 15.7]{I}.

It is natural to formulate Kaplansky's result for set-theoretical complete intersections (trivially, the case $h=0$ gives back Theorem 2, in the local case):

\vspace{.2cm}

\noindent\textbf{Corollary} Let $(R,m)$ be a local noetherian ring, $I,J\subsetneq R$ ideals of $R$, $\underline x=x_1,\ldots,x_h$ a regular sequence in $R$ such that $\sqrt{I\cap J}=\sqrt{\underline xR}$, but $\sqrt I,\sqrt J\neq\sqrt{\underline xR}$.
\[\operatorname{depth}(I+J,R)\leq h+1\]
\hfill$\square$

Though Kaplansky's result can be proved quickly using local cohomology, it is in some sense strong: For example, it implies immediately (by working in the graded instead of in the local situation) that two \lq skew\rq{} (i.\,e., not contained in a plane) lines in $\mathbb{P}^3$ are not a set-theoretical complete intersection, cf. \cite[Example 15.10]{I} and note that this specific example is also covered both by \emph{Faltings' connectedness theorem} (e.\,g. as given in a simple form in \cite[Theorem 15.11]{I}; the full statement \cite[Cor. 4]{F} is more sophisticated) and by \emph{Hartshorne's connectedness theorem} stating that, in projective space over a field, every set-theoretical complete intersection of positive dimension in connected in codimension one (\cite[3.4.6]{H1} or \cite[Theorem 1.3]{H2}; see \cite[Satz 2.6' on p.\,20]{R} for a stronger version).

And though it is strong in the aforementioned sense, Kaplansky's result generalizes smoothly from \emph{set-theoretical} complete intersections to \emph{cohomologically} complete intersections: The sole content of this note is to prove this somewhat surprising fact (theorem in section \ref{section_Result} below).

The proof makes significant use of the fact that, for a cohomologically complete intersection ideal $I$ in a local complete ring $R$, one has
\[\operatorname{End}_R\left(\operatorname{H}^h_I(R)\right)\buildrel\text{canonically}\over=R,\]
where $h$ denotes the codimension of $I$ (\cite[Theorem 2.2\,(iii)]{HeSt}).

\section{Result}

\label{section_Result}An ideal $I$ of a local noetherian ring $R$ is called \textit{cohomologically complete intersection} (\cite{HeSc1}, \cite{HeSc2}) whenever $H^h_I(R)=0$ for all $i\neq h$ for some $h$ (which is then necessarily the height of $I$). In \cite[Theorem 3.2]{HeSc1} it was shown that for $R$ Gorenstein this condition is completely encoded in homological properties of $\operatorname{H}^{\operatorname{height}(I)}_I(R)$. (\cite[Theorem 4.4]{HeSc2} is a related result for modules.)

This condition is clearly weaker than being a \emph{set-theoretical} complete intersection; the actual difference between these two notions can be described in terms of regular sequences on the Matlis dual of the (only) local cohomology module (\cite[1.1.4 Corollary]{He}).

\vspace{.3cm}

\noindent\textbf{Theorem} Let $(R,m)$ be a local noetherian ring, $I,J\subsetneq R$ ideals of $R$ such that $\sqrt I\nsubseteq\sqrt J$, $\sqrt J\nsubseteq\sqrt I$ and $I\cap J$ is a cohomologically complete intersection of depth $h$.
\[\operatorname{depth}(I+J,R)\leq h+1.\]

\textit{Proof:} Assume that $\operatorname{depth}(I+J,R)\geq h+2$. The exact Mayer-Vietoris sequence
\[0\to\LCM_I^h(R)\oplus\LCM_I^h(R)\to\LCM^h_{I\cap J}(R)\to0\]
then leads to a canonical isomorphism
\[\LCM_I^h(R)\oplus\LCM_J^h(R)=\LCM^h_{I\cap J}(R)\]
and hence, by application of the functor $\_\otimes_R\hat R$, to a canonical isomorphism
\[\LCM_{I\hat R}^h(\hat R)\oplus\LCM_{J\hat R}^h(\hat R)=\LCM^h_{(I\cap J)\hat R}(\hat R).\]
Note that $(I\cap J)\hat R$ is a cohomologically complete intersection of depth $h$.

But the endomorphism ring of $\LCM^h_{(I\cap J)\hat R}(\hat R)$ is canonically isomorphic to $\hat R$, by \cite[Theorem 2.2\,(iii)]{HeSt}. Therefore, the direct sum decomposition of $\LCM^h_{(I\cap J)\hat R}(\hat R)$ must be trivial. In particular, either of the two endomorphism rings of $\LCM^h_{I\hat R}(\hat R)$ and $\LCM^h_{I\hat R}(\hat R)$ must be the zero ring, w.\,l.\,o.\,g. we assume this happens for the former one.

This means of course that $\LCM^h_{I\hat R}(\hat R)$ is zero, i.\,e. $\LCM_I^h(R)=0$. But this is only possible if all prime ideals minimal over $I$ have height greater than $h$.

Therefore, no prime ideal minimal over $I$ is minimal over $I\cap J$, i.\,e. to each prime ideal $p$ minimal over $I$ there exists a prime ideal $p_0$ containing $J$ and such that $p\supseteq p_0$, i.\,e. $\sqrt I\supseteq \sqrt J$, contradiction.\hfill$\square$

\end{document}